\documentclass[11pt,a4paper]{amsart}
\usepackage{a4wide}

\usepackage[utf8]{inputenc}
\usepackage[english]{babel}

\usepackage{amsfonts,amssymb,amsmath}

\usepackage{xcolor}
\usepackage[colorlinks,
    linkcolor={red!50!black},
    citecolor={blue!50!black},
    urlcolor={blue!80!black}]{hyperref}

\usepackage{tikz}
\usepackage[all]{xy}
\usetikzlibrary{calc, matrix, arrows, cd}

\newtheorem{theorem}{Theorem}[section]

\theoremstyle{definition}
\newtheorem{definition}[theorem]{Definition}

\newtheorem{example}[theorem]{Example}

\newtheorem{remark}[theorem]{Remark}
\newtheorem*{remark*}{Remark}

\newtheorem*{claim*}{Claim}

\newcommand{\Z}{\mathbb{Z}}
\newcommand{\Q}{\mathbb{Q}}
\newcommand{\R}{\mathbb{R}}
\newcommand{\C}{\mathbb{C}}

\begin{document}
\title{Invariants of $2$-knots}
\begin{abstract}
This short survey, which was written to accompany a minicourse at the BIRS conference ‘Topology in dimension 4.5’,  concerns invariants of knotted $2$-spheres in $S^4$, also known as $2$-knots.
It covers invariants extracted from the algebraic topology of the knot exterior, including Alexander invariants, the Farber-Levine pairing and Casson-Gordon invariants,  as well as gauge theoretic and combinatorial invariants.
Details are scarce and new results inexistant.
\end{abstract}

\author[A.~Conway]{Anthony Conway}
\address{Massachusetts Institute of Technology, Cambridge MA 02139}
\email{anthonyyconway@gmail.com}

\maketitle

\section{Introduction}
\label{sec:Introduction}

The goal of this paper is to survey invariants of $2$-knots, i.e.  of locally flat and smooth embeddings $S^2 \hookrightarrow S^4$.
As there are already several books that thoroughly treat~$2$-knots,  such as~\cite{Hillman2KnotBook, CarterSaito,CarterKamadaSaito, Kamada,KawauchiSurvey},  some justification is needed for the existence of this survey.
\begin{itemize}
\item There is currently renewed interest around the topic of knotted surfaces,  and some new developments are not covered in the aforementioned books. 
\item The references above treat constructions, diagrammatics and invariants,  whereas this article focuses solely on invariants and is scarce on details.
The goal is to describe the state of the art in as few pages as possible.
\end{itemize}

Summarising, this survey is not intended to replace existing books on knotted surfaces (nor could it), merely to complement them and list some recent developments in the field.
Some open questions are also sprinkled throughout the text.

\begin{remark*}
While we focus on embeddings $S^2 \hookrightarrow S^4$, we frequently point out which constructions generalise to surfaces of higher genus and occasionally mention embeddings in other ambient $4$-manifolds.
A couple of other disclaimers are in order.
Firstly, this survey neither treats embeddings of disconnected surfaces nor embeddings of nonorientable surfaces.
We also omit surfaces with nonempty boundary.
Secondly, large portions of the results from Sections~\ref{sec:HomotopyGroups},~\ref{sec:AlexanderInvariants}, and~\ref{sec:CassonGordonInvariants} generalise to other dimensions, but we will avoid systematically pointing this out in the interest of brevity.
Thirdly,  we refrain from making historical comments,   and instead refer the historically minded reader to the previously mentioned books, as well as to~\cite{KervaireWeber}.
Finally,  the reader should keep in mind that listing \emph{all} invariants of $2$-knots is a fool's errand as, after all, the homeomorphism type (and homotopy type) of the complement is an invariant,  as are combinations of preexisting invariants.
\end{remark*}


\subsection*{Organisation}

In Section~\ref{sec:HomotopyGroups},  we consider the homotopy type of $2$-knot exteriors, focusing on the knot group and on higher homotopy groups.
In Section~\ref{sec:AlexanderInvariants}, we review Alexander invariants, the Farber-Levine pairing and  Seifert forms.
In Section~\ref{sec:CassonGordonInvariants}, we briefly discuss Casson-Gordon invariants.
The last two sections concern smooth invariants.
Section~\ref{sec:RochlinGauge} contains a discussion of the Rochlin invariant, $d$-invariants,  Chern-Simons invariants and link Floer cobordism maps, while Section~\ref{sec:QuantumFiniteType} is concerned with combinatorial and quantum invariants, such as quandle cocycle invariants, finite type invariants and invariants coming from Khovanov homology.

\subsection*{Acknowledgements}
These notes were written ahead of the BIRS conference on ‘Topology in dimension 4.5’.  
I am grateful to the organisers for setting up the lectures that led to this survey.
I also wish to thank David Cimasoni, Patrick Orson, Lisa Piccirillo,  Danny Ruberman and Ian Zemke for helpful exchanges.

\subsection*{Conventions}

A \emph{manifold} is understood to mean a compact, connected, oriented, topological manifold.
Embeddings are assumed to be locally flat unless otherwise stated.
A \emph{$2$-knot} is an embedding $S^2 \hookrightarrow S^4$ whereas a \emph{surface knot} refers to an embedding $\Sigma_g \hookrightarrow S^4$ of a genus~$g \geq 0$ surface.
We typically write $K \subset S^4$ for a~$2$-knot and $F \subset S^4$ for a surface knot.
We use~$X_K:=S^4 \setminus \nu(K)$ to denote the exterior of a~$2$-knot, i.e. the exterior of an open tubular neighborhood of $K$ in $S^4$ and, similarly~$X_F$ is reserved for surface knot exteriors.
The closed $4$-manifold obtained by surgery on $K \subset S^4$ is denoted by~$M_K$, i.e. ~$M_K=X_K \cup_\partial (S^1 \times D^3)$.
A \emph{classical knot} refers to a knot in $S^3$.
Given a classical knot $J  \subset S^3$,  we use $\Sigma_n(J)$ to denote the $n$-fold branched of $S^3$ branched along $J$.
When we refer to homotopy groups and covering spaces, it is understood that basepoints have been fixed.

\section{Homotopy groups of the complement}
\label{sec:HomotopyGroups}

The exterior of a $2$-knot $K \subset S^4$ is readily seen to be a homology circle using Alexander duality.
Therefore, just as in classical knot theory,  the homology of the exterior does not provide any useful invariants of $K$ and it makes sense to focus on homotopy groups.

\subsection{The knot group}
\label{sub:KnotGroup}

The \emph{knot group} of a $2$-knot $K \subset S^4$ refers to the group $\pi_1(X_K)$.
A \emph{$2$-knot group} is a group that arises as the knot group of some $2$-knot.
It is known that if $\pi$ is a $2$-knot group, then it is  finitely presented, $H_1(\pi) \cong \Z$, $H_2(\pi)=0$ and $\pi$ is of weight $1$; i.e.  it has an element whose conjugates generate it~\cite{Kervaire}.

\begin{remark}
Kervaire proved that if a finitely presented group $\pi$ satisfies $H_1(\pi) \cong \Z$, $H_2(\pi)=0$ and has deficiency $1$, then it is a $2$-knot group~\cite{Kervaire}.
In contrast with other dimensions,
the deficiency condition is known not to be necessary~\cite{Kervaire}.
\end{remark}

Yajima showed that $2$-knot groups have Wirtinger presentations~\cite{YajimaWirtinger}.
Yajima also showed that the groups of ribbon $2$-knots are characterised by the fact that they have a deficiency one Wirtinger presentation that abelianises to $\Z$~\cite{YajimaRibbon2KnotGroup}.
More tractable group theoretic obstructions to being ribbon can be found in~\cite{CochranRibbon}.
We refer to~\cite{Hillman2KnotBook} for a wealth of information on $2$-knot groups and to~\cite[Section 3.2]{CarterKamadaSaito} for sample calculations.

\begin{remark}
\label{rem:Unknotting}
If a $2$-knot $K  \subset S^4$ has group $\Z$, then it is (topologically) unknotted~\cite{Freedman,FreedmanQuinn}.
In fact knotted surfaces of genus $g$ in $S^4$ with knot group $\Z$ are unknotted for $g \neq 1,2$~\cite{ConwayPowell}; the cases $g=1,2$ remain open.
\end{remark}

Further invariants can be extracted from the knot group.
For example,  the \emph{meridional rank}~\cite{JosephPongtanapaisan} is the minimal number of meridians that generate~$\pi_1(X_K)$, whereas the \emph{Ma-Qiu index}~\cite{MaQiu} is the minimal size of a normal generating set of the commutator subgroup of~$\pi_1(X_K)$,  the \emph{weak unknotting number} (referred to as  the \emph{algebraic stabilisation number} in~\cite{JosephKlugRuppikSchwartz}) is the minimal number of relations of the form~$x=y$, where~$x$ and~$y$ are meridians of~$K$, which abelianise~$\pi_1(X_K)$~\cite{Kanenobu}.
These latter two invariants can be used to give lower bounds on the unknotting number of a $2$-knot~$K$~\cite{HosokawaMaedaSuzuki}.
Here the \emph{unknotting number} of~$K$ (referred to as the \emph{stabilisation number} in~\cite{JosephKlugRuppikSchwartz}) is the minimal number of stabilisations of~$K \subset S^4$ required to obtain an unknotted surface (every $2$-knot $K \subset S^4$ has finite unknotting number~\cite{HosokawaKawauchi}; see~\cite{Boyle,BaykurSunukjian} for related results and~\cite[Appendix A]{ConwayPowell} for the details in the topological category).
The Alexander module~$H_1([\pi_1(X_K),\pi_1(X_K)])=\pi_1(X_K)^{(1)}/\pi_1(X_K)^{(2)}$(and the resulting Alexander invariants) is discussed in Section~\ref{sec:AlexanderInvariants} but we already record that it can used to produce lower bounds on the stabilisation  number~\cite{Miyazaki,MillerPowell,JosephKlugRuppikSchwartz}.

\begin{remark}
\label{rem:Quandle}
The fundamental quandle of a knotted surface is a generalisation of the knot group of a knotted surface.
We refer the interested reader to~\cite{CarterKamadaSaito,FennRourke} for details concerning quandles and racks in the context of knotted surfaces, but simply note for later use that a \emph{quandle} is a set $X$ with a distributive binary operation $*$ that satisfies $a*a=a$ for all $a \in X$ and such that for every $a,b\in X$, there exists a unique $c \in X$ with $c=a*b.$
\end{remark}

\subsection{Higher homotopy groups}
\label{sub:HigherHomotopy}

Contrarily to classical knots whose exteriors are aspherical~\cite{Papakyriakopoulos}, exteriors of $2$-knots can have nontrivial higher homotopy groups; the first published account of this seems to be in~\cite{AndrewsCurtis}.
Gordon later produced pairs of twist spun~$2$-knots~$K,K' \subset S^4$ with isomorphic knot groups but where $\pi_2(X_K)$ and $\pi_2(X_{K'})$ differ as~$\Z[\pi_1]$-modules~\cite{Gordon}.
This was generalised to arbitrarily large families of $2$-knots~\cite{Suciu,PlotnickTheHomotopy} and later to an infinite family~\cite{PlotnickInfinitely}.

\begin{remark}
Epstein proved that if $K$ is a twist spun knot, then $\pi_2(X_K)$ is a free abelian group of infinite rank~\cite{EpsteinLinkingSpheres}.
Andrews and Lomonaco later determined $\pi_2(X_K)$ of a twist spin as a $\Z[\pi_1]$-module~\cite{AndrewsLomonaco,LomonacoSecondProof} and Lomonaco determined $\pi_3(X_K)$,  also for twist spun knots~\cite{Lomonacopi3}.
\end{remark}

Plotnick and Suciu used the $k$-invariant to prove that there are $2$-knots whose exteriors have the same fundamental group,  the same $\pi_2$ and yet are not homotopy equivalent~\cite{PlotnickSuciu}.
Here, the (second) k-invariant of a space $X$ is a cohomology class $k_2(X) \in H^{3}(B\pi_1(X),\pi_2(X))=H^3(\pi_1(X),\pi_2(X))$ that measures how far the second Postnikov stage $P_2(X)$ is from being homotopy equivalent to~$K(\pi_2(X),2) \times B\pi_1(X)$.
The $k$-invariant of a $2$-knot exterior remains challenging to calculate~\cite{LomonacoCompute,PlotnickSuciu}.
The triple $(\pi_1(X),\pi_2(X),k(X))$ is sometimes called the \emph{algebraic $3$-type of $X$.}
Lomonaco proved that for $2$-knots $K \subset S^4$ whose universal cover $\widetilde{X}_K$ satisfies $H_3(\widetilde{X}_K)=0$ (such knots are called \emph{quasi-aspherical}),  the algebraic $3$-type of $X_K$ determines its homotopy type~\cite[Theorem 10.1]{LomonacoCompute}.
Many $2$-knots are quasi-aspherical~\cite[Theorem~10.4]{LomonacoCompute}, but not all~\cite{Ratcliffe}.
Is the homotopy type of the exterior of a~$2$-knot determined by its algebraic $3$-type?
For further open questions concerning the homotopy groups of $2$-knot exteriors, we refer to~\cite[Section XII]{LomonacoCompute}.

 \begin{remark}
 \label{rem:UnknottingSpheres}
Poincar\'e duality (in the universal cover) implies that a $2$-knot $K \subset S^4$ with knot group $\Z$ has aspherical exterior.
The converse was proved in~\cite{DyerVasquez}; see also~\cite{Eckmann} and~\cite[page 31]{Hillman2KnotBook}.
Since $2$-knots $K \subset S^4$ with knot group $\Z$ (or \emph{$\Z$-spheres} for short) are topologically unknotted~\cite{Freedman,FreedmanQuinn},  we conclude that, in the topological category, the unknot is the only $2$-knot with aspherical exterior,  in sharp contrast with classical knot theory.
Is a smooth $2$-knot~$K \subset S^4$ with knot group $\Z$ smoothly unknotted?
The question is also open for surfaces of higher genus; see Remarks~\ref{rem:Unknotting} and~\ref{rem:Exotica} for related considerations.
\end{remark}

\begin{remark}
\label{rem:AsphericalSurgery}
Contrarily to the exterior, the manifold $M_K$ obtained by surgery on a nontrivial~$K \subset S^4$ can be aspherical.
We note that examples are given by fibred $2$-knots with closed fibre an irreducible $3$-manifold $Y$ with infinite $\pi_1$ (such $Y$ are aspherical~\cite{Papakyriakopoulos}) but we refer to~\cite{Hillman2KnotBook,HillmanInventiones} for more details on the homotopy type of $M_K$.\footnote{Given a fixed group $\pi$,  knowing the homotopy type of the surgery on a $2$-knot feeds into applications of (topological) surgery theory to 
classifying $2$-knots with knot group $\pi$ (up to \emph{Gluck reconstruction}); see e.g.~\cite{Hillman2KnotBook,HillmanInventiones}.}
Other examples include those of Cappell-Shaneson which were used to prove that $n$-knots are not determined by their complements for $n=2,3,4$~\cite{CappellShaneson}.
The fact that $2$-knots are not determined by their complement was proved independently by Gordon~\cite{GordonKnotComplement}.
Gluck had previously shown that if a $2$-knot is fixed, then there is at most one other $2$-knot with homeomorphic exterior~\cite{Gluck}.
\end{remark}

The homotopy theoretic invariants considered in this section are also defined for knotted surfaces of higher genus,  but there seems to be less literature on the topic.

\begin{remark}
\label{rem:IntersectionForms}
Given a knotted surface $F \subset S^4$, the $\Z[\pi_1(X_F)]$-module $\pi_2(X_F)=H_2(\widetilde{X}_F)$ is endowed with a sesquilinear Hermitian pairing 
$$H_2(\widetilde{X}_F) \times H_2(\widetilde{X}_F) \to \Z[\pi_1(X_F)]$$
known as the \emph{equivariant intersection form.}
The intersection forms of the various (branched) covers of a knotted surface provide more tractable invariants but, in a simply-connected $4$-manifold $X$, the equivariant intersection form has the advantage of being  a complete invariant of knotted surfaces with knot group $\Z$~\cite{ConwayPowell}; see also~\cite{ConwayPiccirilloPowell}.
In particular,  the equivariant intersection form determines the (oriented) homotopy type of the complement of such \emph{$\Z$-surfaces} and one can again ask about other knot groups.
\end{remark}

\section{Alexander invariants}
\label{sec:AlexanderInvariants}

As in the classical dimension, the homology of the infinite cyclic cover of a $2$-knot exterior provides amenable invariants.
These invariants are often referred to as \emph{Alexander invariants} and include the Alexander module,  the Alexander ideals, the Alexander polynomial, and the Farber-Levine pairing.

\subsection{The Alexander module}
\label{sub:Alexander}

Given a $2$-knot $K \subset S^4$, the kernel of the abelianisation homomorphism $\pi_1(X_K) \twoheadrightarrow H_1(X_K) \cong \Z$ determines an infinite cyclic cover $X_K^\infty \to X_K$.
The deck transformation group of this cover is $\Z$ and therefore the homology groups of $X_K^\infty$ are endowed with the structure of $\Z[t^{\pm 1}]$-modules.

\begin{definition}
The \emph{Alexander module} of a $2$-knot $K \subset S^4$ is the $\Z[t^{\pm 1}]$-module $H_1(X_K^\infty)$.
\end{definition}

The Alexander module is finitely generated,  $\Z[t^{\pm 1}]$-torsion,  and multiplication\footnote{The fact that multiplication by $t-1$ induces an isomorphism $H_1(X_K^\infty) \to H_1(X_K^\infty)$ implies that $H_1(X_K^\infty)$ is $\Z[t^{\pm 1}]$-torsion~\cite[Proposition 1.3]{LevineKnotModules}.} by $t-1$ induces an isomorphism $H_1(X_K^\infty) \to H_1(X_K^\infty)$.
Such $\Z[t^{\pm 1}]$-modules are called \emph{knot modules}~\cite{LevineKnotModules}.
The reason for focusing on the first homology group of $X_K^\infty$ is that $H_0(X_K^\infty) \cong \Z, H_i(X_K^\infty)=0$ for $i \neq 0,1,2$ and $H_2(X_K^\infty)=\operatorname{Ext}^1_{\Z[t^{\pm 1}]}(H_1(X_K^\infty),\Z[t^{\pm 1}])$~\cite[Theorem 2.6]{LevineKnotModules}.
Contrarily to the classical case, the Alexander module of a $2$-knot does not necessarily admit a square presentation matrix.
Also, the Alexander module of a $2$-knot can admit $\Z$-torsion, whereas this is not the case for classical knots.
Levine proved a realisation result for knot modules with no $\Z$-torsion~\cite{LevineKnotModules}.
Ribbon $2$-knots are examples of $2$-knots whose Alexander module has no $\Z$-torsion~\cite{HittExamples}; see~\cite[Theorem 2.2]{CochranRibbon} for a generalisation of this fact.
Kawauchi characterised Alexander modules of ribbon $2$-knots~\cite[Theorem~3.2]{KawauchiTheFirst}.
More information on Alexander modules can be found in~\cite{LevineKnotModules,
GutierrezModules,HillmanAlexanderIdeals, HillmanAlgebraic, KawauchiTheFirst}; this last reference also contains a wealth of information on the Alexander module of knotted surfaces of higher genera.

\begin{example}
\label{ex:AlexanderModule}
For $K \subset S^4$ a fibred $2$-knot with fibre $Y$,  one has~$X_K^\infty \cong Y \times \R$ and so the Alexander module of $K$ satisfies $H_1(X_K^\infty)=H_1(Y)$.
For example,  for $n \neq 0,n$-twist-spinning a classical knot~$J \subset S^3$ results in a fibred $2$-knot~$K=\tau^n(J) \subset S^4$ with closed fibre~$\Sigma_n(J)$~\cite{Zeeman} and Alexander module~$H_1(X_K^\infty)=H_1(\Sigma_n(J))$.
\end{example}

\begin{definition}
\label{def:AlexanderPolynomial}
The \emph{Alexander polynomial} $\Delta_K \in \Z[t^{\pm 1}]$ of a $2$-knot $K$ refers to the order of~$H_1(X_K^\infty)$; it is well defined up to multiplication by $\pm t^k$ for $k \in \Z$.
\end{definition}

We refer to~\cite{TuraevIntroductionTo} for the definition of the order of a module over a Noetherian factorial domain,  but in order to give a sense of the idea, we recall its definition for modules over~$\Q[t^{\pm 1}]$.
Since $\Q[t^{\pm 1}]$ is a PID,  any torsion $\Q[t^{\pm 1}]$-module $H$ splits as a direct sum of cyclic modules~$\Q[t^{\pm 1}]/(p_i(t))$.
The \emph{order} of the~$\Q[t^{\pm 1}]$-module $H$ is the product of the $p_i(t)$.


As for classical knots, the Alexander polynomial of a $2$-knot $K$ satisfies $\Delta_K(1) =\pm 1$~\cite{LevineKnotModules} but contrarily to the classical case,  it need not be symmetric.
In fact Kinoshita proved that every polynomial $f(t) \in \Z[t^{\pm 1}]$ with $f(1)=\pm 1$ arises as the Alexander polynomial of a ribbon $2$-knot~\cite{Kinoshita}; see also~\cite{YajimaRibbon2KnotGroup,HillmanRibbonLink} for shorter proofs.
Nakanishi and Nishizawa gave a topological condition on a 2-knot ensuring that its Alexander polynomial is symmetric~\cite{NakanishiNishizawa}.
Moussard and Wagner provide conditions for the Alexander polynomial of certain ribbon $2$-knots to factor as $f(t)f(t^{-1})$~\cite{MoussardWagner}.

\begin{remark}
Further invariants can be extracted from the Alexander module.
The \emph{$k$-th elementary Alexander ideal} of a $2$-knot refers to the $k$-th elementary ideal of its Alexander module; here we refer the reader to~\cite{TuraevIntroductionTo} for the definition of elementary ideals and note that conventions surrounding these ideals may vary.
The $0$-th ideal (also known as the \emph{Alexander ideal}) of a classical knot $J \subset S^3$ is principal and is generated by the Alexander polynomial.
The Alexander ideal of a $2$-knot need not be principal; this was recently exploited by Joseph in the study of $0$-concordance of $2$-knots~\cite{Joseph}.
The \emph{$k$-th Alexander polynomial} $\Delta_K^{(k)}$ of a $2$-knot~$K$ refers to the greatest common divisor of its $k$-th Alexander ideal; we refer to~\cite{LevineKnotModules} for more properties of these polynomials.
Further Alexander invariants include the \emph{Nakanishi index} of a $2$-knot (the minimal number of generators of~$H_1(X_K^\infty)$~\cite{Nakanishi,JosephKlugRuppikSchwartz}) and the \emph{determinant} of a $2$-knot (the evaluation of the Alexander ideal of $K$ at $t=-1$~\cite{Joseph}).
All of these Alexander invariants can be generalised to surface knots of higher genus.
\end{remark}

\begin{remark}
\label{rem:Normalisation}
The Alexander polynomial of a link is only defined up to multiplication by~$\pm t^n$ with~$n \in \Z$.
In the classical dimension,  if one wishes to work with a symmetric representative of $\Delta_L(t)$, then the indeterminacy further reduces to fixing a sign.
Fixing this sign leads to the definition of Conway's potential function and its resulting skein theory~\cite{ConwayEnumeration}.
Such a normalisation is currently not known to exist for the Alexander polynomial of $2$-links.
Even for $2$-knots, there is no canonical normalisation  because $\Delta_K$ is generally not symmetric (in~\cite[Section~2]{HabiroKanenobuShima} a choice of a normalisation is fixed for ribbon $2$-knots, but it is not canonical).
Can an analogue of Conway's potential function be defined for arbitrary $2$-links?
We refer to Giller's work~\cite[Section~5]{Giller} for some considerations related to this question that involve Seifert matrices (the definition of these matrices is recalled in Subsection~\ref{sub:Seifert}).
Giller focuses on $2$-knots whose (generic) projections to $\R^3$ have no triple points (i.e.  ribbon $2$-knots~\cite{YajimaSimply,KawauchiPseudo}).
\end{remark}

\begin{remark}
\label{rem:TwistedAlexander}
In the classical dimension,  twisted Alexander polynomials have proved to be successful generalisations of the Alexander polynomial; see~\cite{FriedlVidussi} for a  survey.
The definition of these invariants extends to $n$-knots and for $n=2$, they have been studied by Kanenobu and Sumi with applications to the enumeration of ribbon $2$-knots~\cite{KanenobuSumiClassification2018,
KanenobuSumiClassification2019,KanenobuSumiClassification2020}.
\end{remark}


%

\subsection{The Farber-Levine form} 
\label{sub:FarberLevine}

Given a $\Z[t^{\pm 1}]$-module $H$, use $T_\Z H$ to denote the $\Z$-torsion subgroup of $H$.
The \emph{Farber-Levine pairing} of a $2$-knot~$K \subset S^4$, introduced independently by Farber~\cite{FarberRussian,FarberDuality} and Levine~\cite{LevineKnotModules},  is a non-singular symmetric~$\Z$-bilinear form
$$ \operatorname{FL}_K \colon T_\Z H_1(X_K^\infty) \times T_\Z H_{1}(X_K^\infty) \to \Q/\Z. $$
The Farber-Levine pairing satisfies $FL_K(tx,ty)=FL_K(x,y)$ for every $x,y \in  T_\Z H_1(X_K^\infty)$.
\begin{remark}
Here is a slightly informal definition of $FL_K([\alpha],[\beta])$.
Since $T_\Z H_1(X_K^\infty)$ is finite,  there exists an integer $m \in \Z$ and a surface $F \subset X_K^\infty$ with $\partial F=m \alpha$.
Levine shows that~$(t^k-1)H_2(X_K^\infty;\Z_m)=0$ for $k>0$ large enough,  so that $(t^k-1)F=\partial Y +m F_0$ for some $3$-manifold $Y \subset X_K^\infty$ and some surface $F_0 \subset X_K^\infty$~\cite{LevineKnotModules}.
Now $FL_K([\alpha],[\beta]):=\frac{-1}{m}(Y \cdot \beta)$, where~$\cdot$ denotes the algebraic intersection number.
\end{remark}

The Farber-Levine pairing can be calculated using the linking form of any Seifert solid for~$K$~\cite[Theorem 7.3]{FarberDuality}; see also~\cite[Section 4]{Ogasa}.
Here, the \emph{linking form} of a $3$-manifold $Y$ refers to the pairing
$$\ell_Y \colon TH_1(Y) \times TH_1(Y) \to \Q/\Z,$$
where $TH_1(Y)$ denotes the torsion subgroup of $H_1(Y)$ and, roughly speaking, $\ell_Y([a],[b])=\frac{1}{n}(F \cdot b)$, where $F \subset Y$ is any surface with boundary the disjoint union of $n$ copies of $a$.

\begin{example}
\label{ex:FarberLevine}
For a fibred $2$-knot $K \subset S^4$,  the Farber-Levine pairing $FL_K$ corresponds to the linking form on the fibre~\cite[Proposition 7.1]{LevineKnotModules}; see also~\cite{FarberRussian,FarberDuality}.
In particular,  for $n \neq 0$, the Farber-Levine pairing of the $n$-twist spin on a classical knot $J \subset S^3$ is $\ell_{\Sigma_n(J)}$.
\end{example}

The Farber-Levine pairing was applied by Farber to exhibit a $3$-knot group that is not a~$2$-knot group~\cite{FarberRussian}.
The Farber-Levine pairing has also been applied to the study of double sliceness~\cite{StoltzfusDouble,StoltzfusIsometries},  as well as to obstruct invertibility~\cite{FarberRussian, FarberDuality} and ribbon move equivalence~\cite{Ogasa}.
The Farber-Levine pairing can be extended to higher genus surfaces~\cite[The Second Duality Theorem]{KawauchiThreeDualities}; see also~\cite{Sekine,KawauchiIntrinsic,
KawauchiLinkingSignature,KawauchiPseudo} for further applications, generalisations and invariants derived from this pairing.

\subsection{Seifert matrices}
\label{sub:Seifert}

In classical knot theory,  Seifert matrices lead to presentation matrices for the Alexander module and, in particular,  can be used to calculate the Alexander polynomial.
The usefulness of Seifert matrices in $2$-knot theory is more restricted as we now describe.
Given an abelian group $A$,  we write $TA \subset A$ for its torsion subgroup and~$FA=A/TA$ for its torsion-free part.
The \emph{$\pm$-Seifert form} of an embedded $3$-manifold $Y \subset S^4$ is defined as 
$$\beta_\pm \colon H_2(Y) \times FH_1(Y) \ \to \Z, ([F],[\gamma]) \mapsto \ell k(F,\gamma^{\pm })$$
where $\gamma^\pm \subset S^4$ denotes a curve obtained by pushing $\gamma$ off $Y$ in the $\pm$-direction.
The fact that the $\pm$-Seifert form is well defined follows from Alexander duality.
Seifert forms are not knot invariants, altough they provide insights into the Alexander module, as we now describe.

\emph{Seifert matrices} $V^+$ and $V^-$ for a $2$-knot $K \subset S^4$ are matrices for the $\pm$-Seifert form of a Seifert solid $Y$ of  $K$ with respect to bases of $FH_1(Y)$ and~$H_2(Y)$.
Note that~$H_2(Y)$ is free,  but $H_1(Y)$ need not be free.
The $\pm$-Seifert forms extend over~$\Q$,  leading to matrices~$V^\pm_\Q$.
The matrix~$tV^+_\Q-V^-_\Q$ presents the rational Alexander module~$H_1(X_K^\infty;\Q)$~\cite{LevineCodimensionTwo}.
If $H_1(Y)$ is torsion free, then $tV^+-V^-$ presents the Alexander module $H_1(X_K^\infty)$; see e.g.~\cite{Giller,MoussardWagner}. In particular, in this case,   the Alexander module admits a square presentation matrix.

\begin{remark}
\label{rem:GromovNorm}
While not directly related to Seifert matrices, we note that Ruberman defined the \emph{Gromov norm} $|K|$ of a $2$-knot $K \subset S^4$ as the infimum of the Gromov norm of all Seifert solids for $K$~\cite{RubermanSeifert}.
Here the Gromov norm of a real singular chain $z \in \sum r_i\sigma_i$ in a space $X$ is~$|z|:=\sum_i |r_i|$ and the Gromov norm of a manifold $X$ is the Gromov norm of its fundamental class. 
Among other applications,  Ruberman proves that $|K|=0$ if $K$ is a ribbon $2$-knot.
\end{remark}

\section{Casson-Gordon invariants}
\label{sec:CassonGordonInvariants}

We recall the definition of the Casson-Gordon $3$-manifold invariants~\cite{CassonGordon1,CassonGordon2} and outline how they give rise to $2$-knot invariants.
Given a closed $3$-manifold $Y$ and a homomorphism $\chi \colon H_1(Y) \to \Z_d$, bordism theory implies that there exists a pair $(W,\psi)$, where $W$ is a $4$-manifold with boundary~$n$ disjoint copies of $Y$,  and $\psi \colon H_1(W) \to \Z_d$ is a homomorphism extending~$\chi$.
The \emph{Casson-Gordon invariant} of $(Y,\chi)$ is defined as 
$$\sigma(Y,\chi):=\frac{1}{n}(\sigma^\psi(W)-\sigma(W)) \in \Q.$$
Here $\sigma^\psi(W)$ denotes the signature of the twisted intersection form on $H_2(W;\C^\psi)$.
Note that~$\sigma(Y,\chi)$ is a special case of the Atiyah-Singer $U(1)$-invariant~\cite{AtiyahSingerEllipticIII}.

Recall that $M_K$ denotes the closed $4$-manifold obtained by surgery on a $2$-knot~$K \subset S^4$.
Use $\overline{Y} \subset M_K$ to denote the result of capping off a Seifert solid $Y \subset S^4$ for $K$.
Pick a lift of~$\overline{Y}$ to the infinite cyclic cover $M_K^\infty$ of  $M_K$.
We again denote this lift by $\overline{Y}$ and write~$\iota \colon H_1(\overline{Y}) \to H_1(M_K^\infty)$ for the inclusion induced map.
Given a prime power order homomorphism $\chi \colon H_1(M_K^\infty) \to \Z_d$, Ruberman shows that~$ \sigma(K,\chi):=\sigma(\overline{Y},\chi \circ \iota) \in \Q$
does not depend on the choice of the Seifert solid~$Y$; see~\cite[Theorems 4.7 and 4.8]{RubermanDoublySlice} and~\cite[Theorem~2.1]{RubermanCG}.

\begin{definition}
\label{def:CassonGordon}
Let $d$ be a prime power.
The \emph{Casson-Gordon invariant} of a $2$-knot $K \subset S^4$ associated to a homomorphism $\chi \colon H_1(M_K^\infty)\to \Z_d$ is defined as 
$$ \sigma(K,\chi):=\sigma(\overline{Y},\chi \circ \iota) \in \Q,$$
where $\overline{Y} \subset M_K$ is any capped off Seifert solid for $K$.
\end{definition}

\begin{example}
\label{ex:CG}
For $K \subset S^4$ a fibred $2$-knot,  one has~$M_K^\infty \cong \overline{Y} \times \R$,  with $\overline{Y}$ a capped-off Seifert solid for $K$,  and therefore~$\sigma(K,\chi)=\sigma(\overline{Y},\chi)$.
In particular for $K=\tau^n(J)$ the $n$-twist spin on a classical knot $J \subset S^3$, with $n \neq 0$, one has $\sigma(K,\chi)=\sigma(\Sigma_n(J),\chi)$.
Casson-Gordon $3$-manifold invariants can be calculated using the surgery formulas from~\cite[Lemma 3.1]{CassonGordon1},~\cite[Theorem 6.7]{CimasoniFlorens}, and~\cite[Theorem 3.6]{GilmerConfiguration}.
\end{example}

In~\cite{RubermanDoublySlice} Ruberman gives a different (but equivalent) definition of these invariants using infinite cyclic covers.
He also shows that if $K$ is ribbon, then $\sigma(K,\chi)=0$ for all $\chi$.
This relies on the fact that every ribbon $2$-knot bounds a punctured $\#_{i=1}^n S^1 \times S^2$ for some $n$~\cite{Yanagawa}; see also~\cite[Theorem 5.2]{RubermanDoublySlice}.
Ruberman also uses Casson-Gordon invariants to obstruct double sliceness,  amphicheirality and invertibility of $2$-knots, while Ogasa has applied these invariants to obstruct ribbon move equivalence~\cite{Ogasa}.

\begin{remark}
Given a $4$-manifold $X$ with the integral homology of $S^1 \times S^3$, Ruberman has recently introduced a Casson-Gordon type invariant for knotted tori $T \subset X$ with $H_1(T) \to H_1(X)$ surjective~\cite{RubermanTori}.
Assuming that $X$ is smooth, Ma~\cite{Ma} later related this invariant to the Furuta-Ohta invariant~\cite{FurutaOhta} of $X$ and to a gauge theoretic invariant of smooth knotted tori defined by Echeverria~\cite{Echeverria}.
\end{remark}

\section{The Rochlin invariant and gauge theoretic invariants}
\label{sec:RochlinGauge}

We describe how the Rochlin invariant, the $d$-invariant from Heegaard-Floer homology and the Chern-Simons function, all of which are $3$-manifold invariants give rise to invariants of $2$-knots.
We conclude with a brief description of link Floer cobordism maps.
Up to now all invariants were defined for locally flat surfaces, whereas the invariants in this section require smooth embeddings.

\subsection{The Rochlin invariant}
\label{sub:Rochlin}

We recall the definition of the Rochlin invariant and outline how it gives rise to a smooth $2$-knot invariant.
Given a $3$-manifold $Y$ and a spin structure~$\mathfrak{s}$ on $Y$,  there exists a smooth spin $4$-manifold $W$ and a spin structure $\mathfrak{s}_W$ on $W$ such that~$\partial (W,\mathfrak{s}_W)=(Y,\mathfrak{s})$ (i.e.~$\Omega_3^{\operatorname{Spin}}=0$).
A theorem of Rochlin~\cite{Rochlin} ensures that the signature of a closed smooth spin $4$-manifold is divisible by $16$.
We can therefore define the \emph{Rochlin invariant} of $(Y,\mathfrak{s})$ as 
$$ \mu(Y,\mathfrak{s}):=\sigma(W) \in \Z_{16}.$$
We describe how $\mu$ defines a smooth $2$-knot invariant.
Given a Seifert solid $Y \subset S^4$ for a smooth 2-knot $K$,  use $\mathfrak{s}$ to denote the spin structure on $\overline{Y} \subset M_K$ obtained by restricting the spin structure from $S^4$ to $Y$ and then extending it over $\overline{Y}$.
Ruberman~\cite{RubermanDoublySlice} showed that the Rochlin invariant $\mu(\overline{Y},\mathfrak{s})$ only depends on $K$; see also~\cite[Theorem 7]{Sunukjian0Concordance}.
\begin{definition}
\label{def:Rochlin}
The \emph{Rochlin invariant} of a smooth $2$-knot $K \subset S^4$ is defined as 
$$\mu(K)=\mu(\overline{Y},\mathfrak{s})\in \Z_{16}$$
where $\overline{Y}$ is any capped off Seifert solid for $K$.
\end{definition}

\begin{example}
\label{ex:Rochlin}
The Rochlin invariant of a fibred $2$-knot is the Rochlin invariant of its closed fibre with respect to the aforementioned spin structure $\mathfrak{s}$.
For example,  the closed fibre of~$K=\tau^2(T_{3,5})$, the $2$-twist spin of the $(3,5)$-torus knot, is the branched cover $\Sigma_2(T_{3,5})$ also known as the Poincar\'e homology sphere~\cite{KirbyScharlemann}.
It follows that $H_1(X_K^\infty)=0$ and $FL_K=0$ whereas~$\mu(K)=1$.
Here we used that integer homology spheres have a single spin structure.
\end{example}

The Rochlin invariant vanishes on ribbon 2-knots~\cite[Theorem 5.3]{RubermanDoublySlice}, obstructs amphicheirality and double sliceness~\cite{RubermanCG, RubermanDoublySlice},  as well as $0$-concordance~\cite{Sunukjian0Concordance}, and ribbon move equivalence~\cite{Ogasa}.
Can the Rochlin invariant be defined for nullhomologous $2$-knots
 in other closed simply-connected smooth $4$-manifolds?
 If so,  can it distinguish topologically isotopic $2$-knots?

\subsection{$d$-invariants}
\label{sub:dInvariants}

Correction terms (also known as $d$-invariants) are rational numbers that are extracted from the Heegaard-Floer homology package~\cite{HeegaardFloer}.
Traditionally,  a rational number $d(Y,\mathfrak{s}) \in \Q$ is associated to a rational homology $3$-sphere~$Y$ and a spin-c structure $\mathfrak{s}$ on~$Y$~\cite{OzsvathSzaboAbsolutely}.
This was extended to $3$-manifolds with $b_1(Y)>0$ using twisted coefficients~\cite{BehrensGolla,LevineRubermanCodimension1}.
The details are beyond the scope of the survey apart from the fact, due to Levine and Ruberman, that if~$X$ is a smooth $4$-manifold with the integral homology of~$S^1 \times S^3$, $Y \subset X$ is a $3$-dimensional submanifold,
and~$\mathfrak{s}_X$ is the spin-c structure on $Y$ induced from $X$, then a correction term~$d(Y,\mathfrak{s}_X) \in \Q$ can be defined and only depends on $X$ and on the homology class~$y:=[Y] \in H_3(X)$~\cite[Theorem~1.1]{LevineRubermanCodimension1}; the resulting invariant is denoted $d(X,y)$.
Taking~$X$ to be~$M_K$, the result of surgery on a smooth $2$-knot~$K \subset S^4$, leads to the following definition, which is due to Levine and Ruberman~\cite{LevineRubermanCodimension1}.

\begin{definition}
The \emph{$d$-invariant} $d(K) \in \Q$ of a smooth $2$-knot $K  \subset S^4$ is the
 $d$-invariant $d(M_K,y) \in \Q$, where $y\in H_3(M_K)=\Z$ is the homology class that is Poincar\'e dual to the class in~$H^1(M_K)=\operatorname{Hom}(H_1(M_K),\Z)$ that maps the meridian of $K$ to $1 \in \Z$. 
\end{definition}

\begin{example}
\label{ex:d}
The $d$-invariant of a fibred $2$-knot $K \subset S^4$ can be calculated in terms of (twisted) $d$-invariants of a closed fibre of $K$.
For example, the $2$-twist spin of the $(3,5)$-torus knot,  denoted $K=\tau^2(T_{3,5})$, has closed fibre the Poincar\'e homology sphere (recall Example~\ref{ex:Rochlin}) and therefore $d(K)=2$~\cite[Example 5.1]{LevineRubermanCodimension1}.
\end{example}
The $d$-invariant vanishes for ribbon $2$-knots  and obstructs invertibility and amphicheirality~\cite[Section 6]{LevineRubermanCodimension1}; it also obstructs $0$-concordance~\cite{Sunukjian0Concordance}.
Do other homology cobordism invariants of $3$-manifolds lead to invariants of $2$-knots?

\subsection{The Chern-Simons invariant}
\label{sub:ChernSimon}

Given a closed~$3$-manifold~$Y$, the \emph{Chern-Simons functional} is a function~$CS_Y \colon \mathcal{A}/\mathcal{G} \to \R/\Z$, where~$\mathcal{A}=\Omega^1(Y;\mathfrak{su}_2)$ denotes the set of connections on a trivialised principal $SU(2)$-bundle~$P=Y \times SU(2)$ and~$\mathcal{G}=C^\infty(Y,SU(2))$ is the gauge group.
Here is one  possible definition of $CS_Y$.
Any connection~$A \in \mathcal{A}$ extends to a connection~$\widetilde{A}$ on some smooth~$4$-manifold~$W$ with boundary $Y$ and, using~$F_{\widetilde{A}}$ to denote the curvature of $\widetilde{A}$,  the \emph{Chern-Simons functional} of $A$ is defined as
$$ CS_Y(A):=\frac{1}{8\pi^2}\int_W \operatorname{Tr}\left( F_{\widetilde{A}} \wedge F_{\widetilde{A}}\right).$$
Given a closed smooth~$4$-manifold~$X$ with~$H_*(X)=H_*(S^1 \times S^3)$ and a nonzero cohomology class~$c \in H^1(X)$,  Taniguchi constructs a Chern-Simons function~$cs_{X,c} \colon \widetilde{R}(X) \to \R/\Z$ on the set~$\widetilde{R}(X)$ of flat connections on~$P_X=X \times SU(2)$ considered up to  degree $0$ gauge equivalence~\cite{TaniguchiInstantons}.
%
The value of~$CS_{X,c}$ on a flat connection~$a \in \widetilde{R}(X)$ is
$$ CS_{X,c}(a):=-\frac{1}{8\pi^2}\int_{\widetilde{X}^c} \operatorname{Tr}(F_{A_a} \wedge F_{A_a}),$$
where $p \colon \widetilde{X}^c \to X$ is the~$\Z$-cover determined by~$c \in  H^1(X)$, and~$A_a$ is a connection on~$P_{\widetilde{X}^c}=\widetilde{X}^c \times SU(2)$ that, informally speaking, agrees with~$p^*a$ on the left of~$\widetilde{X}^c$ and vanishes on the right of~$\widetilde{X}^c$;
see~\cite{TaniguchiInstantons,TaniguchiSeifert} for details.
Note that  while $\widetilde{X}_c$ is noncompact,  the integral~$CS_{X,c}(a)$ is finite because it can be expressed in terms of $CS_Y$ for a closed $3$-dimensional submanifold~$Y \subset X$ that is Poincar\'e dual to~$c$~\cite[Lemma 4.2]{TaniguchiInstantons}.

The composition~$H_1(X) \xrightarrow{c} \Z \to \Z_j$  determines a $\Z_j$-cover~$p_j \colon X_{j,c} \to X$, and therefore a functional~$cs^j_{X,c}:=cs_{X_{j,c},p_j^*(c)}$, where~$p_j^*(c) \in H^1(X_{j,c})$.
Taking~$X$ to be~$M_K$, the result of surgery on a smooth~$2$-knot~$K \subset S^4$, leads to the following definition, which is due to Taniguchi~\cite{TaniguchiSeifert}.
\begin{definition}
\label{def:ChernSimon}
For $j \in \Z$, the~\emph{$j$-th Chern-Simons invariant} of a smooth~$2$-knot~$K$ is the function~$cs^j_{M_K,c} \colon \widetilde{R}((M_K)_{j,c}) \to \R/\Z=(0,1]$,  where $c \in H^1(M_K)=\operatorname{Hom}(H_1(M_K),\Z)$ is the class that maps the meridian of $K$ to $1 \in \Z$. 
\end{definition}

Taniguchi proves that for every $j \in \Z$, the set $\operatorname{im}(cs^j_{M_K,c})$ is finite and $\operatorname{im}(cs^j_{M_K,c})=\lbrace 1\rbrace$ if~$K$ is ribbon~\cite{TaniguchiSeifert}.
These invariants also obstruct invertibility~\cite[Theorem 1.1]{TaniguchiSeifert}, and
sample calculations 
can be found in~\cite[Proposition~1.3 and Example 1.4]{TaniguchiSeifert}.

\begin{remark}
Other gauge theoretic invariants give rise to smooth $2$-knot invariants.
For example,  Mrowka, Ruberman and Saveliev's~\cite{MrowkaRubermanSaveliev} invariant $\lambda_{SW}(X)$ of smooth homology $S^1 \times S^3$'s gives a smooth $2$-knot invariant by setting $\lambda_{SW}(K):=\lambda_{SW}(M_K)$.
We will not list all $2$-knot invariants that can be obtained this way: any closed $4$-manifold invariant gives a $2$-knot invariant via the surgery manifold.
\end{remark}

\begin{remark}
\label{rem:Exotica}
Gauge theory has proved fruitful at distinguishing topologically isotopic surfaces that are smoothly embedded 
in a closed
smooth $4$-manifold.
In order to briefly survey the state of the art on this topic, we call two smoothly knotted surfaces \emph{exotic} if they are topologically but not smoothly isotopic.
The first\footnote{If one drops the requirement that the surfaces be oriented,  then the first examples are due to Finashin, Kreck and Viro who constructed exotic embeddings of $\#_{i=1}^{10} \R P^2$ into $ S^4$~\cite{FinashinKreckViro}.} examples of exotic surfaces are due to Fintushel-Stern~\cite{FintushelSternRimSurgery,
FintushelSternAddendum} 
who constructed exotic symplectically embedded surfaces $F_0,F_1 \subset X$ in a simply-connected symplectic $4$-manifold~$X$ where $\pi_1(X \setminus F_i)=1$; for other knot groups, see~\cite{KimTwistSpinning,
KimRubermanSmooth}.
The surfaces are constructed by rim surgery on a surface $\Sigma \subset X$ and are distinguished using the relative Seiberg-Witten invariants $SW_{X,\Sigma}^{\mathcal{T}}$; see also~\cite{Mark} for the use of Heegaard-Floer theory.
Other constructions of exotic surfaces rely on the existence of closed exotic $4$-manifolds~\cite{HoffmanSunukjian, TorresSmoothly,
TorresTopologically,BenyahiaMalechTorres}.
Similarly to the unknotting conjecture in $S^4$ (recall Remark~\ref{rem:UnknottingSpheres}),  one can ask:
is there a smoothly nontrivial nullhomologous 
$2$-knot $K\subset \C P^2$
 that is topologically unknotted?
It is also unknown whether there are exotic $2$-knots whose homology class in $H_2(\C P^2) \cong \Z$ is $d=1,2$.
For $d>3$,  these classes cannot be represented by smoothly embedded spheres~\cite{KronheimerMrowka} but their minimal genus representatives admit exotic copies~\cite{KimTwistSpinning,Finashin}.
\end{remark}

\subsection{Cobordism maps in link Floer homology}

We outline very briefly how the functoriality of link Floer homology gives rise to invariants of closed surfaces in $S^4$~\cite{ZemkeSetUp,ZemkeFunctoriality}.
For earlier work on cobordism maps in link Floer homology, we refer to~\cite{Juhasz,JuhaszMarengon}.
Given the technically demanding nature of these papers, we only describe the theory in very loose terms in the case of knots,
but note that it extends to links.
Associated to a doubly based knot~$\mathbb{K}$ in a $3$-manifold $Y$,  there is a doubly filtered \emph{curved chain complex}~$\mathcal{CFL}^\infty(Y,\mathbb{K})$ over the ring~$\mathbb{F}_2[U,V]$~\cite{ZemkeSetUp}.
Here $\mathbb{F}_2$ denotes the field with two elements.
Given doubly based knots~$\mathbb{K}_1 \subset  Y_1$ and $\mathbb{K}_2 \subset Y_2$,  Zemke shows that a decorated smooth knot cobordism~$\mathcal{F}=(\Sigma,\mathcal{A})$ in a smooth cobordism~$(W^4;Y_1,Y_2)$
and a spin-c structure~$\mathfrak{s}$ on~$W$ induce a~$\Z \oplus \Z$-filtered homomorphism of~$\mathbb{F}_2[U, V]$-modules
$$ F_{W,\mathcal{F},\mathfrak{s}} \colon \mathcal{CFL}^\infty(Y_1,\mathbb{K}_1,\mathfrak{s}|_{Y_1}) \to \mathcal{CFL}^\infty(Y_2,\mathbb{K}_2,\mathfrak{s}|_{Y_2})$$
which is an invariant up to~$\mathbb{F}[U,V]$-equivariant,~$\Z \oplus \Z$-filtered chain homotopy.
A decorated smoothly embedded surface~$F \subset S^4$ can be viewed as a smooth cobordism from the empty link in $S^3$ to itself and therefore gives a homomorphism
$$ F_{S^4,\mathcal{F},\mathfrak{s}} \colon \mathcal{CFL}^\infty(\emptyset) \to \mathcal{CFL}^\infty(\emptyset).$$
There is a canonical identification~$\mathcal{CFL}^\infty(\emptyset) \cong \mathbb{F}_2[U^{\pm 1},V^{\pm 1}]$ and therefore $ F_{S^4,\mathcal{F},\mathfrak{s}}$ is determined by $F_{S^4,\mathcal{F},\mathfrak{s}}(1)$.
Unfortunately this map uniquely  depends on the genus of the surface~\cite[Theorem 1.8]{ZemkeFunctoriality}.
Knot Floer homology can be used to smoothly distinguish topologically isotopic surfaces with non-empty boundary~\cite{JuhaszMillerZemke,DaiMallickStoffregen}.

\section{Quantum and finite type invariants}
\label{sec:QuantumFiniteType}

In this final section, we survey invariants that arise from quantum and combinatorial means.

\subsection{Quandle cocycle invariants}
\label{sub:CocycleInvariants}

Associated to a quandle~$X$ (recall Remark~\ref{rem:Quandle} for the definition) and an abelian group~$A$, there are quandle cohomology groups~$H_Q^*(X;A)=Z_Q^*(X;A)/B_Q^*(X;A)$.
The definition is reminiscent of group cohomology, but we refer to~\cite{CarterJelsovskyKamadaLangfordSaito} for details.
Given a diagram of a smoothly embedded surface~$F \subset S^4$ (whose set of sheets is denoted~$\mathcal{R}$) and a finite quandle~$X$, a \emph{colouring} refers to a function~$\mathcal{C} \colon \mathcal{R} \to X$ which satisfies certain conditions at the double point set~\cite[Section 4.1.3]{CarterKamadaSaito}.
Given a~$3$-cocycle~$\theta \in Z^3_Q(X,A)$, Carter, Jelsovsky, Kamada, Langford and Saito~\cite{CarterJelsovskyKamadaLangfordSaito} introduced the \emph{quandle cocycle invariant}
\begin{equation}
\label{eq:CocycleInvariant}
\Phi_\theta(F):= \sum_{\mathcal{C}}\prod_\tau B(\tau,\mathcal{C}) \in \Z[A].
\end{equation}
Here the product is taken over all triple points~$\tau$ of the diagram,  the sum is taken over all possible colourings~$\mathcal{C}$, and the \emph{Boltzmann weight}~$B(\tau,\mathcal{C}):=\theta(x,y,z)^{\varepsilon(\tau)} \in A$ is defined combinatorially in terms of the colouring of the diagram around the triple point~$\tau$; see~\cite{CarterJelsovskyKamadaLangfordSaito,CarterKamadaSaito} for details.
Implicit in~\eqref{eq:CocycleInvariant} is the fact that~$\Phi_\theta(F)$ is independent of the choice of the diagram for~$F$~\cite{CarterJelsovskyKamadaLangfordSaito}; this is proved using that any two diagrams for~$F$ are related by a sequence of Rosemann moves~\cite{RosemanReidemeister} (these are analogues of the Reidemeister moves in the context of diagrams of knotted surfaces).

The value of~$\Phi_\theta(F)$ only depends on the cohomology class~$[\theta] \in H_Q^3(X;A)$~\cite{CarterJelsovskyKamadaLangfordSaito}.
The cocyle invariant has been applied to reprove that certain twist-spun knots are not invertible~\cite{CarterJelsovskyKamadaLangfordSaito}; (see also~\cite{SatohSurfaceDiagrams,
AsamiSatoh,CarterElhamdadiGranaSaito}) as well as to study triple point numbers~\cite{SatohShima}.
A variant of the quandle cocycle invariant can be used to obstruct ribbon concordance~\cite{CarterSaitoSatohRibbon}; see~\cite{NosakaQuandleHomotopy} for another variant.
Can one define quandle cocycle invariants for properly embedded surfaces with nonempty boundary?
Can quandle cocycle invariants smoothly distinguish topologically isotopic knotted surfaces, e.g. the exotics discs from~\cite{Hayden}?


\subsection{Finite type invariants of 2-knots}
\label{sub:FiniteType}

Building on the concept of finite type invariants in classical knot theory~\cite{Vassiliev,BirmanLin,Bar-Natan},
Habiro, Kanenobu and Shima defined finite type invariants for smooth ribbon~$2$-knots~\cite{HabiroKanenobuShima} and showed that the coefficients of the normalised Alexander polynomial (recall Remark~\ref{rem:Normalisation}) are of finite type.
In later work, Habiro and Shima determined the rational vector space of rational-valued finite type invariants and proved that it only depends on the coefficients of the normalised Alexander polynomial~\cite{HabiroShima}.


\subsection{Khovanov-Jacobsson numbers}
\label{sub:Khovanov}

Khovanov homology associates a graded abelian group $\operatorname{Kh}(L)$ to a classical link~$L \subset S^3$~\cite{Khovanov}.
It is known to be functorial, meaning that for every smooth cobordism~$C \subset S^3 \times [0,1]$ between links~$L_0,L_1 \subset S^3$ there is a homomorphism~$\Phi_C \colon \operatorname{Kh}(L_0) \to \operatorname{Kh}(L_1)$~\cite{Khovanov, Jacobsson,ClarkMorrisonWalker}.
A smooth knotted surface~$F \subset S^4$ can be viewed as a cobordism from the empty link to itself and therefore gives a homomorphism~$\phi_F \colon \operatorname{Kh}(\emptyset) \to \operatorname{Kh}(\emptyset)$, i.e.  an endomorphism of~$\Z$, i.e. a integer.
The integer~$\Phi_F$ is called the \emph{Khovanov-Jacobsson number} of the knotted surface~$F \subset S^4$.
Unfortunately, ~$\Phi_F$ was shown to equal~$2^{g(F)}$ where~$g(F)$ denotes the genus of~$F$~\cite{Tanaka}; see also~\cite{Rasmussen} as well as~\cite{CarterSaitoSatoh} for related results.
Thus Khovanov-Jacobsson numbers cannot distinguish knotted surfaces of the same genus (does a similar statement hold for Khovanov-Rozansky $\mathfrak{gl}(n)$-homology~\cite{KhovanovRozansky}?)
This is in stark contrast with the case of surfaces with boundary,  where recent work of Hayden and Sundberg shows that certain relative Khovanov-Jacobsson numbers can smoothly distinguish topologically isotopic discs~\cite{HaydenSundberg}; see also prior work of Sundberg and Swann on Khovanov-Jacobsson \emph{classes}~\cite{SundbergSwann}.

\bibliographystyle{plain}
\bibliography{2knotInvariants}
\end{document}